\documentclass[11pt]{imsart}

\usepackage{amssymb, amsmath, amsthm, srcltx} 
\usepackage{graphicx} 
\usepackage{color} 
\usepackage{bbm}

\startlocaldefs
\newtheorem{thm}{Theorem}[section] 
\newtheorem{lem}[thm]{Lemma} 
\usepackage{pst-plot} 
\usepackage{pst-grad} 
\newtheorem{cor}[thm]{Corollary} 
\newtheorem{prop}[thm]{Proposition} 
\theoremstyle{definition} 
\newtheorem{defn}[thm]{Definition} 
\newtheorem{rem}[thm]{Remark} 
\numberwithin{equation}{section} 
\newcommand{\E}{\mathbf{E}\,} 
\renewcommand{\P}{\mathbf{P}} 
\newcommand{\G}{\mathcal{G}} 
\newcommand{\pS}{\mathcal{S}} 
\newcommand{\ver}{\mathcal{V}}

\newcommand{\PS}{\mathcal{P}} 
\newcommand{\bi}{\mathbf{i}}
 
\newcommand{\con}[1]{Cont(#1)} 
\newcommand{\Xb}{\mathbf{X}}
\newcommand{\Tr}{\operatorname{Tr}}
\newcommand{\iii}{i_0,\dots,i_{2mp-1}} 
\newcommand{\bigprod}{\prod_{j=0}^{2mp-1}{X^{\varepsilon}_{i_j i_{j+1}}}} 
\endlocaldefs

\begin{document}

\begin{frontmatter}

\title{Asymptotic distribution of singular values of powers of random matrices}
\runtitle{Powers of random matrices}

\author{\fnms{Nikita} \snm{Alexeev}}\sep 
\author{\fnms{Friedrich} \snm{G\"otze}}
\and
\author{\fnms{Alexander} \snm{Tikhomirov}\thanksref{rf}}
\thanksref{crc}
\thankstext{rf}{Partially supported by RF grant of the leading scientific schools
NSh-638.2008.1. Partially supported RFBR, grant  N 09-01-12180 and RFBR--DFG, grant N 09-01-91331.}
\thankstext{crc}{Partially supported by CRC 701 ``Spectral Structures and Topological
Methods in Mathematics'', Bielefeld}
\runauthor{N. Alexeev\sep F. G\"otze \and A. Tikhomirov}

\affiliation{Faculty of Mathematics and Mechanics, Saint-Petersburg State University, Russia}
\affiliation{Faculty of Mathematics, University of Bielefeld, Germany} 
\affiliation{Department of Mathematics Komi Research Center of Ural Branch of RAS, Syktyvkar State University, Russia} 

  \begin{abstract}
Let $x$ be a complex random variable such that ${\E {x}=0}$, ${\E |x|^2=1}$, ${\E |x|^{4} < \infty}$.
Let $x_{ij}$, $i,j \in \{1,2,\dots\}$ be independet copies of $x$.  Let ${\Xb=(N^{-1/2}x_{ij})}$, $1\leq i,j \leq N$ be a random matrix. Writing $\Xb^*$ for the adjoint matrix of $\Xb$, consider the product $\Xb^m{\Xb^*}^m$
with some $m \in \{1,2,\dots \}$. The matrix $\Xb^m{\Xb^*}^m$ is Hermitian positive semi-definite. Let $\lambda_1,\lambda_2,\dots,\lambda_N$ be eigenvalues of $\Xb^m{\Xb^*}^m$ (or squared singular values of the matrix $\Xb^m$ ). 
In this paper we find the asymptotic distribution function 
\[
 G^{(m)}(x)=\lim_{N\to\infty}\E{F_N^{(m)}(x)}
\]
of the empirical distribution function 
\[
{F_N^{(m)}(x)} = N^{-1} \sum_{k=1}^N {\mathbb{I}{\{\lambda_k \leq x\}}},
\]
where $\mathbb{I} \{ A\}$ stands for the indicator function of event $A$. The moments of $G^{(m)}$ satisfy
\[
M^{(m)}_p=\int_{\mathbb{R}}{x^p\, dG^{(m)}(x)}=\frac{1}{mp+1}\binom{mp+p}{p}.
\]
In Free Probability Theory $M^{(m)}_p$ are known as Fuss--Catalan numbers.
With $m=1$ our result turns to a well known result of Marchenko--Pastur 1967.
\end{abstract} 


\end{frontmatter}
\section{Introduction} 
Let $\Xb=(N^{-1/2}x^{(N)}_{ij}), 1\leq i,j \leq N$ be a random matrix. We assume that  $x_{ij}\equiv x_{ij}^{(N)}$ are independent complex random variables such that 
\begin{equation}
\E {x_{ij}}=0,\quad \E |x_{ij}|^2=1,\quad \E |x_{ij}|^{4}\leq B 
\label{cond_mom}
\end{equation}
 with some $B<\infty$ independent of $N$. We assume additionally that 
\begin{equation}
L_N(\alpha)=N^{-2}\sum_{1\leq i,j\leq N}\E|x_{ij}|^4\mathbb{I}\{|x_{ij}|>\alpha\sqrt{N} \} \to 0\textrm{ as }N\to\infty
\label{cond_ratio}
\end{equation}
for all $\alpha>0$. Note that $x_{ij}\equiv x^{(N)}_{ij}$ and $\Xb \equiv \Xb^{(N)}$ can depend on $N$, which is not reflected in our further notation.

Writing $\Xb^*$ for the adjoint matrix of $\Xb$, consider the product
\[
\mathbf{W}^{(m)}=\Xb^m{\Xb^*}^m
 \]
with some $m \in \{1,2,\dots \}$. The matrix $\mathbf{W}^{(m)}$ is Hermitian positive semi-definite. Let $\lambda_1,\lambda_2,\dots,\lambda_N$ be eigenvalues of $\mathbf{W}^{(m)}$ (or squared singular values of the matrix $\Xb^m$ ). 
In this paper we find the asymptotic distribution function 
\[
 G^{(m)}(x)=\lim_{N\to\infty}\E{F_N^{(m)}(x)}
\]
of the empirical distribution function 
\[
{F_N^{(m)}(x)} = N^{-1} \sum_{k=1}^N {\mathbb{I}{\{\lambda_k \leq x\}}},
\]
where $\mathbb{I} \{ A\}$ stands for the indicator function of event $A$. 
\begin{thm}
Assume that $(\ref{cond_mom})$ and $(\ref{cond_ratio})$ hold. Then the limit $G^{(m)}(x)=\lim_{N\to\infty}\E{F_N^{(m)}(x)}$ exists. The function $G^{(m)}(x)$ is a distribution function and it has moments 
\begin{equation}
M^{(m)}_p=\int_{\mathbb{R}}{x^p\, dG^{(m)}(x)}=\frac{1}{mp+1}\binom{pm+p}{p}.
\label{moments}
\end{equation}
\label{main_thm}
\end{thm}

\begin{cor}
Let $x_{ij}$ be independet copies of a random variable, say $x$, such that 
\[
\E {x}=0,\quad \E |x|^2=1,\quad \E |x|^{4}< \infty.
\]
Let $\Xb=(N^{-1/2}x_{ij})$, $1\leq i,j \leq N$. 
Then the limit ${\lim_{N\to \infty}\E F_N^{(m)}(x)}$ exists and it is equal to $G^{(m)}(x)$.
\label{iid} 
\end{cor}
Gessel and Xin 2006 \cite{Ges2006} showed that for any natural $m$ the sequence $M^{(m)}_1$, $M^{(m)}_2, \dots$ is a sequence of moments of some probability measure. Hence,  $G^{(m)}$ is a probability distribution for any natural $m$.  Since $M^{(m)}_p \leq c_m^p$    with some $c_m<\infty$, by Carleman's Theorem  in \cite{Carl1926} the measure $G^{(m)}$ is uniquely determined by its moments.  The support of the measure $G^{(m)}$ is the interval $\left[0,m^{-m}(m+1)^{m+1}\right]$. 
%

With ${m=1}$ Theorem \ref{main_thm} turns to a well known result of Marchenko--Pastur 1967~\cite{MarPas}. Namely, the asymptotic distribution $G^{(1)}$ of eigenvalues of the matrices  $\Xb\Xb^*$ has the moments ${M^{(1)}_p=\frac{1}{p+1}\binom{2p}{p}}$. Note that in the case $m=1$ our fourth moment assumption
is stronger than assumptions in Theorem 2.5 and Theorem 2.8 in Bai 1999 \cite{Bai99}. The question of the weakest sufficient conditions in the case $m>1$ remains an open problem.

In Free Probability Theory $M^{(m)}_p$ are known as Fuss--Catalan numbers. Combinatorial properties of this sequence have been studied by Nica and Speicher~2006~\cite{NSp}.  Mlotkowski 2009~\cite{Ml09} investigated a family of distributions, say $G^{(m,r)}$, with real $m\geq 0$ and ${0\leq r \leq m}$, such that  $G^{(m,r)}$ has moments $\frac{r}{mp+p+r}\binom{mp+p+r}{p}$. It is easy to check that ${G^{(m,1)}=G^{(m)}}$. 
Oravecz 2001 \cite{ora} proved that powers of Voiculescu's circular element have distribution~$G^{(m)}$. This distribution belongs to the class of Free Bessel Laws  (see Banica~et~al~2008~\cite{banica08}). 

Let $\mathcal{M}_m(x)=\sum_{p=0}^\infty{M^{(m)}_p x^p}$ be the generating function of the sequence $M^{(m)}_p$.  
It satisfies the following functional equation (see equation (7.68) on p. 347  Graham et al 1988 \cite{CMath})
\begin{equation} 
\mathcal{M}_m(x)=1+x \mathcal{M}_m^{m+1}(x). 
\label{fe} 
\end{equation} 
Equation (\ref{fe}) allows us to  describe $G^{(m)}$ in 
 the framework of Free Probability Theory. In Free Probability Theory the  free multiplicative convolution $\xi \boxtimes \eta$ is defined for any positive random variables $\xi$ and $\eta$ (see Nica and Speicher 2006 \cite{NSp}, p. 287). The $S$-transform is a
homomorphism with respect to free multiplicative convolution, i.e. if $\xi$ and $\eta$ are free independent positive variables, then $S_{\xi \boxtimes \eta}(z)=S_{\xi}(z)S_\eta(z)$. 
Recall that the $S$-transform, say $S(z)$, of a distribution $\mu$ is defined as follows. Let
\begin{equation} 
M_p=\int_{\mathbb{R}}{x^p\, d \mu(x)}, \quad u(z)=\sum_{p=1}^{\infty} M_p z^p.
\nonumber
\end{equation}
Then
\begin{equation} 
S(z)=\frac{z+1}{z}u^{-1}(z),  
\end{equation} 
where $u^{-1}$ denotes  the inverse function of $u$.

Equation (\ref{fe}) allows to calculate the $S$-transform, say $S^{(m)}(z)$, of $G^{(m)}$, and 
\begin{equation} 
S^{(m)}(z)=\frac{1}{(1+z)^m}. 
\end{equation} 
It means that the family $G^{(m)}$ has 
the following  property: if a random variable  $\xi$ has distribution $G^{(m)}$ then the $r$-th power of the $S$-transform of  $\xi$ is equal
 to the $S$-transform of multiplicative free power $\xi^{\boxtimes r}$. 
This property holds for this family of distributions only. 
 

To prove Theorem \ref{main_thm} we use truncation and the method of moments. Truncation means that we can replace (see Section \ref{sec_trunc} for details) $\Xb$ by the matrix ${\widetilde{\Xb}=(\widetilde{X}_{ij})}$ with truncated entries (here and below ${X_{ij}=N^{-1/2}x_{ij}}$ denote entries of matrix $\Xb$)
\begin{equation}
\widetilde{X}_{ij}=X_{ij}\mathbb{I}\{|X_{ij}|<\alpha_N \},
\end{equation}
where $\alpha_N$ is some sequence of positive numbers such that $\alpha_N \to 0$ as $N \to \infty$. 
Lemma~\ref{trunc_lem} (see Section \ref{sec_trunc}) reduces the proof of Theorem~\ref{main_thm} to the proof of the next proposition. 
\begin{prop}
\label{pr_trunc}
Assume that $\alpha_N \to 0$ and $\beta_N \to 0$. Then Theorem $\ref{main_thm}$ holds if
\begin{equation}
\left|X^{(N)}_{ij}\right|\leq \alpha_N, \quad \max_{1\leq i,j \leq N} \left|\E X^{(N)}_{ij}\right|\leq \beta_N N^{-3/2}, \quad \left|\E \left|X^{(N)}_{ij}\right|^2-1/N \right| \leq \beta_N N^{-3/2}.
\label{cond}
\end{equation}
\end{prop}
Let us explain our proof of Proposition \ref{pr_trunc}. Denote by $\xi_{m}(N)$ a random variable with distribution $\E F_{N}^{(m)}$. We show that the moments $\E \xi^p_{m}(N)$ converge to $M^{(m)}_p$. 
In order to simplify the notation assume for a while that ${X}_{ij}$ are real random variables.
 Then one can represent $\E \xi^p_{m}(N)$ as 
\begin{equation}
\E \xi^p_{m}(N)={\sum}^{(2mp)} N^{-1} {\E \prod_{j=0}^{2mp-1} X^{\varepsilon(j)}_{i_j i_{j+1}}},
\nonumber
\end{equation}
where the sum ${\sum}^{(2mp)}$ is taken over ${{i_0,..,i_{2mp} \in \left\lbrace 1,..,N \right\rbrace}}$ such that $i_{2mp}=i_0$.
The notation $X^{\varepsilon(j)}_{i_j i_{j+1}}$ means $X^+_{i_j i_{j+1}} := X_{i_{j} i_{j+1}}$ in case of $\varepsilon(j)=+$ and $X^-_{i_j i_{j+1}} := X_{i_{j+1} i_{j}}$ in case of $\varepsilon(j)=-$ (see Section \ref{sec_mom} for a precise definition of the spin variable $\varepsilon(j)$). We investigate properties of paths  ${(i_0,..,i_{2mp})}$ by combinatorial methods. The moment $\E \xi^p_{m}(N)$ converges to the number of paths of a special type. Namely, one can describe such paths as follows: the cardinality of $\{i_0,..,i_{2mp}\}$ is equal to $mp+1$ and each factor $X_{i_j i_{j+1}}$ appears in the product $\bigprod$ twice. In Section \ref{sect_count} we count the number of these paths.
\section{The proof of the main result.} 
\subsection{Truncation.} 
\label{sec_trunc}
%
Recalling that $X_{ij}=N^{-1/2}x_{ij}$, we can rewrite $L_N(\alpha)$ as
\[
 L_N(\alpha) = \sum_{1\leq i,j\leq N}\E|X_{i,j}|^4\mathbb{I}\{|X_{i,j}|>\alpha \}.
\]
Since for all $\alpha>0$ the ratio $L_N(\alpha)/\alpha^{4}$ tends to $0$, one can find a sequence $\alpha_N \downarrow 0$ such that $L_N(\alpha_N)/\alpha_N^{4} \to 0$ and $N^\delta \alpha_N^{-1}\to \infty$ for any $\delta>0$ as $N \to \infty$.
Let $\widetilde{F}^{(m)}_N(t)$ denote the empirical spectral distribution function of the matrix ${\widetilde{\Xb}}^m{{{\widetilde{\Xb}}}}^{*m}$. 
\begin{lem} 
The limit behaviors of  $\E\widetilde{F}^{(m)}_N(t)$ and $\E{F}^{(m)}_N(t)$ are the same, that is
\[
 \sup_{t \in\mathbb{R}}|\E \widetilde{F}^{(m)}_N(t)-\E F^{(m)}_N(t)|\to 0.
\]

\label{trunc_lem}
\end{lem}
\begin{proof} 
Since by definition  $|\widetilde{F}^{(m)}_N(t)-F^{(m)}_N(t)| \ne 0$ only if there exist $i,j \in \{1,\dots,N\}$ such that $|X_{ij}|\geq\alpha_N$,
 we have 
\begin{equation} 
|\E \widetilde{F}^{(m)}_N(t)-\E F^{(m)}_N(t)|\leq \sum_{1\leq i,j\leq N}{\P(|X_{ij}|\geq \alpha_N)}. 
\label{tr_est} 
\end{equation} 
Estimating ${\P(|X_{ij}|\geq \alpha_N)\leq \alpha_N^{-4} \E |X_{ij}|^{4} \mathbb{I}\{|X_{ij}|>\alpha_N \}}$ and using inequality (\ref{tr_est}) we obtain 

 \begin{eqnarray}
\nonumber
 \sup_{t \in \mathbb{R}}|\E\widetilde{F}^{(m)}_N(t)-\E F^{(m)}_N(t)| &\leq& \alpha_N^{-4}\sum_{i,j=1}^N\E |X_{ij}|^{4}\mathbb{I}\{|X_{ij}|>\alpha_N \}\\
 &=&\alpha_N^{-4} L_N(\alpha_N) \to 0.  \\ \label{lemtr}
\nonumber 
\end{eqnarray}
\end{proof} 
 
Note, that the lower order moments of the truncated variables are asymptotically equal to the moments of the original variables. Writing for a while  $X=X_{ij}$ we have for $k\leq 3$ 
\begin{equation} 
|\E {\widetilde{X}}^k - \E X^k| \leq \E |X|^k \mathbb{I}\{|X|>\alpha_N\}. 
\label{inlow} 
\end{equation} 
The right hand side of (\ref{inlow}) can be estimated as  
\begin{equation} 
 \E |X|^k \mathbb{I}\{|X|>\alpha_N\} \leq \alpha_N^{k-4}\E |X|^{4} \leq \beta_N N^{-3/2}, 
\end{equation} 
where $\beta_N = B\alpha_N^{k-4}N^{-1/2}\to 0$ as $N \to \infty$.
 
Lemma \ref{trunc_lem} shows that the limit behaviors of $\widetilde{F}^{(m)}_N(t)$ and
${F}^{(m)}_N(t)$ are the same. Thus we may replace $\Xb$ by $\widetilde{\Xb}$ in the following arguments and assume that $\Xb$ is truncated, that is, that entries of $\Xb$ satisfy the assumption (\ref{cond}).


\subsection{Moments of the spectral distribution.}
\label{sec_mom}
We apply the method of moments. Recall that
 ${\lambda_1,\lambda_2,\dots,\lambda_N}$ denote the eigenvalues of $\Xb^m{\Xb^*}^m$. We can write
\begin{equation} 
 \E \xi^p_m(N)=N^{-1}\E {\sum_{j=1}^N {\lambda_j^p}}=N^{-1}\E {\Tr(\Xb^m {\Xb^*}^m)^p}  \; . 
\label{1} 
\end{equation} 

We assume that $m$ and $p$ are fixed and study  the asymptotics  of $\E \xi^p_m(N)$ as $N \rightarrow \infty$. In order to simplify notation, hence forth we assume that $X_{ij}$ are real random variables.

In the  Hermitian case, the trace of $\Xb^{2k}$ 
may be rewritten  in terms of the entries of $\Xb$  via
\begin{equation} 
 \E \Tr \Xb^{2k}= {\sum}^{(2k)} \E{\prod_{j=0}^{2k-1}{{X}_{i_j i_{j+1}}}},
\end{equation} 
where the sum ${\sum}^{(s)}$ is taken over ${{i_0,..,i_{s} \in \left\lbrace 1,..,N \right\rbrace}}$ such that $i_{s}=i_0$. 

In the non-Hermitian case $\E \Tr(\Xb^m{\Xb^*}^m)^p$ has a similar
representation. An entry of $\Xb^m\Xb^{*m}$ is given by
\begin{eqnarray}
[\Xb^m\Xb^{*m}]_{ik}&=&\sum_{1\leq i_j \leq N}X_{ii_1}X_{i_1i_2}\cdots X_{i_{m-1}i_m}X_{i_{m+1}i_m}\cdots X_{k i_{2m-1}}\\
\label{linal}
\end{eqnarray}
We write  $X^+_{i_j i_{j+1}} := X_{i_{j} i_{j+1}}$ and $X^-_{i_j i_{j+1}} := X_{i_{j+1} i_{j}}$. Then the right hand side of (\ref{linal}) takes the form
\begin{equation}
[\Xb^m\Xb^{*m}]_{ik}=\sum_{1\leq i_j \leq N}\prod_{j=0}^{2m-1}X^{\varepsilon(j)}_{i_j i_{j+1}}, 
\end{equation}
where $i_0=i$, $i_{2m}=k$, and the 'spin' variable $\varepsilon(j)$ takes values
 $\varepsilon(j)=+$ with $j<m$, and $\varepsilon(j)=-$ with $j\geq m$.
 Since $(\Xb^m\Xb^{*m})^p=\Xb^m\Xb^{*m}\cdots\Xb^m\Xb^{*m}$ ($p$ times), one needs to change the order of indices in $X^\varepsilon_{i_ji_{j+1}}$ if the spin $\varepsilon=-$ and 
\begin{equation}
\varepsilon(j) = 
\begin{cases} 
 + \text{ , if } j\pmod{2m} \in \{ 0,\dots\,m-1 \}, \\ 
 - \text{ , if } j\pmod{2m} \in \{ m,\dots\,2m-1 \}.
\end{cases}
\label{spin} 
\end{equation}


Using these notions (\ref{1}) takes  the form 
\begin{eqnarray}
\nonumber \E \xi^p_{m}(N)&=&N^{-1} \E \Tr(\Xb^m {\Xb^*}^m)^p\\
 &=&{\sum}^{(2mp)} N^{-1} {\E \prod_{j=0}^{2mp-1} X^{\varepsilon(j)}_{i_j i_{j+1}}}. \label{trace} \\
\nonumber
\end{eqnarray}
A crucial notion in the proof is that of  'paths' of indices of the type $(i_0,i_1,\dots,i_{2mp-1})$.
 
\subsection{Description of  paths.} 
We consider a path $\mathbf{i}=(\iii)$ which corresponds to a product ${\prod_{j=0}^{2mp-1} X^{\varepsilon}_{i_j i_{j+1}}}$. Let $\PS$ be a set of pairs $\{(j,j+1)^{\varepsilon(j)}\}_{j=0}^{2m-2}\bigcup\{(2mp-1,0)^-\}$,  where ${(j,j+1)^{+}:=(j,j+1)}$, $(j,j+1)^{-}:=(j+1,j)$ and $\varepsilon(j)$ is given by (\ref{spin}). We call pairs $(j,j+1)^{\varepsilon(j)}$ and $(k,k+1)^{\varepsilon(k)}$ equivalent (denoted by~$(j,j+1)^{\varepsilon(j)}\sim(k,k+1)^{\varepsilon(k)}$) iff ${X^{\varepsilon(j)}_{i_j i_{j+1}}\equiv X^{\varepsilon(k)}_{i_k i_{k+1}}}$. We also call $(j,j+1)^{\varepsilon(j)}$ an edge of the path $\bi$. 
We construct a directed graph $\G_\bi$ as follows. A vertex $\ver$ of $\G_\bi$ is a subset of ${\{0,1,\dots, 2mp-1\}}$ such that $j\in\ver$ and $k \in \ver$ if and only if $i_j=i_k$. There exists an edge $(\ver,\mathcal{U})$ if and only if  there exist $l\in \ver$ and $r\in \mathcal{U}$ such that $(l,r)\in\PS$ (note that $|l-r|=1$). Denote by $V$ the total number of vertices of the graph~$\G_\bi$ and by $E$ its total number of edges. Since the graph $\G$ is connected $E\geq V-1$. It is clear that $V$ is a cardinality of $\{i_0,i_1,\dots,i_{2m-1}\}$ and $E$ is a cardinality of a quotient set~$\PS/\sim$. Denote by $k_r$ ($r=1,\dots,E$) the cardinality of each equivalence class in $\PS$. Note, that $k_1+k_2+\dots+k_E=2mp$.

\begin{rem} 
Consider paths $\bi=(\iii)$ and $\mathbf{k}=(k_0,k_1,\dots,k_{2m-1})$ such that $\G_\bi=\G_\mathbf{k}$. It is clear that if $x_{ij}$ are identically distributed then 
\[
\E \prod_{j=0}^{2mp-1} X^{\varepsilon(j)}_{i_j i_{j+1}}=\E \prod_{j=0}^{2mp-1} X^{\varepsilon(j)}_{k_j k_{j+1}}.
\]
 We will show, that assuming  our conditions the asymptotic products corresponding to equivalent paths
are equal as well.
\end{rem}

\begin{defn}
We define the contribution of a graph $\G$ to $\ref{trace}$ as
\[
\con{\G}=\sum_{\bi :\G_\bi=\G} {N^{-1} \E \prod_{j=0}^{2mp-1} X^{\varepsilon(j)}_{i_j i_{j+1}}}
\]
\end{defn}


\begin{lem} 
Using these notations we have that the contribution of the path $\G$ is asymptotically given by
\begin{equation} 
\con{\G} \sim 
N^{V-1} \prod_{r=1}^{E}{\E x_{i_s i_t}^{k_r}},
\label{cont} 
\end{equation} 
when $N$ tends to infinity.
\end{lem} 
\begin{proof} 
Since $X_{ij}$ are independent we have 
\begin{equation}
\E \bigprod=\prod_{r=1}^{E}{\E X_{i_s i_t}^{k_r}}.
\nonumber 
\end{equation}
Furthermore, for any vertex $\ver$ the number of possible values of corresponding indices (indices $i_j$ such that $j\in\ver$) lies  between $N$ and ${N-2mp \sim N}$.
 The lower bound $N-2mp$ is due to the fact that indices corresponding to this vertex 
should not coincide with indices corresponding to other vertices and that there are
at most  $2mp$ different indices. This yields the multiplicity $N^{V}$. Together with 
the factor  ${N^{-1}}$ this finally leads to the formula (\ref{cont}). 
\end{proof} 
 
\begin{defn} 
We call  a graph $\G_\bi$ $(m,p)$-regular graph, if it
has at least $mp+1$ vertices and $k_r\geq 2$ for all $r \in \{1,2,\dots,E\}$.
The path $\bi$ we call $(m,p)$-regular path.
\end{defn}

\begin{lem} 
 $\con{\G_\bi}$ does not converge to zero if and only if $\G_\bi$ is the regular path. 
\label{contr_zero} 
\end{lem} 
\begin{proof} 
Since the variables $X_{ij}$ satisfy conditions (\ref{cond}), we have
\begin{equation} 
\left|\E X^k_{i j}\right| \leq \E X^2_{i j}|X_{ij}|^{k-2} \leq N^{-1} \alpha_N^{k-2}.
\label{est_prod} 
\end{equation} 
Of course, this estimaton holds for $k=1$ too.
At first we consider that one of $k_r$ is equal to $1$ (without loss of generality $k_1=1$). Then we have 
\begin{eqnarray}
\nonumber \left|\prod_{r=1}^{E}{\E X_{i(r)j(r)}^{k_r}}\right|&=&\left|\E X_{i(1)j(1)} \prod_{r=2}^{E}{\E X_{i(r)j(r)}^{k_r}}\right|\leq \beta_N N^{-3/2}N^{-E+1}\alpha_N^{\sum_r{(k_r-2)}}\\ 
&=&\beta_N N^{-3/2} N^{-E+1}\alpha_N^{2mp-1-2(E-1)}\leq N^{-E-1/2},\\
\nonumber
\end{eqnarray}

and the contribution of such a graph is bounded by
\begin{equation} 
|\con{\iii}|\leq N^{V-1}N^{-E-1/2}=N^{V-E-1}N^{-1/2}.
\label{in1}
\end{equation}
Note that $V-E-1\leq0$ since the graph $\G$ is connected and hence $N^{V-E-1}N^{-1/2}$ tends to $0$.

Furthermore, we consider the case $V<mp+1$.  Note that  $k_r\geq 2$ for any $r$ and $E \leq {2mp}/{2}=mp$.
Our truncation leads to 
\begin{equation}
\left|\prod_{r=1}^{E}{\E X_{i(r)j(r)}^{k_r}}\right| \leq N^{-E}\alpha_N^{\sum_r{(k_r-2)}}=N^{-E}\alpha_N^{2mp-2E}.
\label{inee} 
\end{equation} 
Using inequality (\ref{inee}) to estimate the terms in (\ref{cont}), we obtain for such a product
\begin{equation} 
N^{V-1} \left|\prod_{r=1}^{E}{\E X_{i(r)j(r)}^{k_r}}\right| \leq  N^{V-E-1}\alpha_N^{2mp-2e}.
\label{contto0}
\end{equation}

Note that $E\geq V-1$  and $2mp-2E\geq 0$. It follows that the right hand side of (\ref{contto0}) does not converge to $0$ only if $2mp-2E=0$ and $V-E-1=0$, i.e. $V=mp+1$ and the graph $\G$ is a regular graph.
\end{proof}
 Furthermore, we obtain 
\begin{lem}
A regular graph is a tree and it has exactly $V=mp+1$ vertices and exactly 
$E=mp$ edges $($ each representing an equivalence class of size~$k_r=2 )$.
\end{lem} 

\begin{rem} 
Due to the fact that $\E X_{ij}^2\sim 1/N$ and by the remarks above we can write the contribution of a regular graph $\G_{reg}$ as

\begin{equation} 
\con{\G_{reg}} \sim 1. 
\label{cont_reg} 
\end{equation}  
 
\end{rem} 

We now show the connection between the moments of the spectral distribution $\E F^{(N)}_m$ 
and the number of regular graphs. Indeed, $\xi_m(N)$ has  distribution $\E F^{(N)}_m$. 
Denote by $T_{m,p}$ the set of all possible graphs of view $\G_\bi$ and by $T^{reg}_{m,p}$ 
the set of all $(m,p)$-regular graphs. Then 
\begin{equation} 
\E \xi^p_m(N)=\sum_{\pS \in T_{m,p}}\con{\pS} \sim \sum_{\pS \in T^{reg}_{m,p}}1=\#T^{reg}_{m,p}. 
\label{eq_reg}
\end{equation}  
We can reformulate \ref{eq_reg} as
\begin{lem} 
$\lim_{N \to \infty}{\E \xi_{N}^p}$ is equal to the number of $(m,p)$-regular graphs. 
\label{lem1} 
\end{lem} 

 \subsection{Counting of the number of regular graphs.}
\label{sect_count}
 \begin{lem} 
  The number of all $(m,p)$-regular graphs is  
 $\# T^{reg}_{m,p}=M^{(m)}_p$.
 \label{lem2} 
 \end{lem} 
 \begin{proof} 
  The numbers $M^{(m)}_p=\frac{1}{p+1}\binom{mp+p}{p}$ satisfy to the recurrence (see \cite{CMath}): 
 \begin{equation} 
 M^{(m)}_p)=\sum_{p-1}{\prod_{i=0}^{m-1}M^{(m)}_{p_i}},~~~M^{(m)}_1=1,
 \end{equation} 
where the sum $\sum_{p-1}$ is taken over all ${p_0+p_1+\dots+p_{m}=p-1}$.
 We will show that there is one-to-one correspondence between collections of $(m,p_k)$-regular graphs $(\G_{m,p_0},\dots,\G_{m,p_{m}}):
\sum_{i=0}^{m} p_i=p-1$ and $(m,p)$-regular graphs $\mathcal{G}_{m,p}$. It follows that 
 \begin{equation} 
  \#T^{reg}_{m,p}=\#\bigcup_{p-1}{\bigotimes_{i:=0}^{m}T^{reg}_{m,p_i}}= 
 \sum_{p-1}{\prod_{i:=0}^{m}{\#T^{reg}_{m,p_i}}} 
 \end{equation} 
 and the sequence $\#T^{reg}_{m,p}$ satisfies to both the same reccurence and initial conditions as the sequence Fuss--Catalan numbers $M^{(m)}_p$ and, by this reason, these two sequences are equal.
 \begin{prop} 
 The number $\# T^{reg}_{m,1}=1$ for all $m$. If $\G_\bi$ is a $(m,1)$-regular graph then indices $i_k$ and $i_l$ are equal iff $(k+l)=2m$.
 \label{3} 
 \end{prop} 
 \begin{proof} 
 By induction. Consider $m=1$. In this case it is clear, that there is only one regular graph $0 \rightarrow 1$ and Proposition \ref{3} holds. 
Assume, that Proposition \ref{3} holds for all $m < m_0$. Consider the path $\bi$ and a corresponding graph $\G_\bi$ (see fig.\ref{pict3}). 
 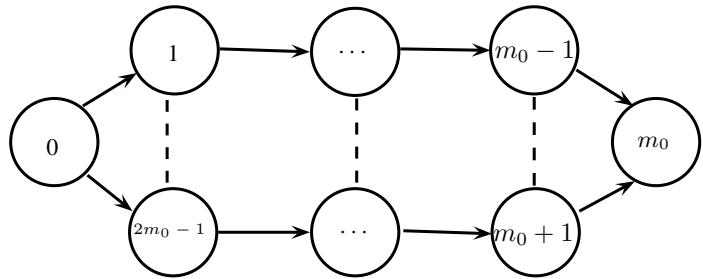
\begin{figure}[hbt] 
   \begin{center} 
      \scalebox{1} 
  { 
  \begin{pspicture}(0,-1.81)(9.2,1.81) 
  \pscircle[linewidth=0.04,dimen=outer](0.6,-0.01){0.6} 
  \pscircle[linewidth=0.04,dimen=outer](2.2,1.21){0.6} 
  \pscircle[linewidth=0.04,dimen=outer](4.6,1.19){0.6} 
  \pscircle[linewidth=0.04,dimen=outer](6.98,1.21){0.6} 
  \pscircle[linewidth=0.04,dimen=outer](8.6,-0.01){0.6} 
  \pscircle[linewidth=0.04,dimen=outer](2.18,-1.21){0.6} 
  \pscircle[linewidth=0.04,dimen=outer](4.6,-1.21){0.6} 
  \pscircle[linewidth=0.04,dimen=outer](7.0,-1.21){0.6} 
  \psline[linewidth=0.04cm,arrowsize=0.1cm 2.0,arrowlength=1.4,arrowinset=0.4]{->}(0.98,0.47)(1.7,0.91) 
  \psline[linewidth=0.04cm,arrowsize=0.1cm 2.0,arrowlength=1.4,arrowinset=0.4]{->}(2.8,1.23)(4.0,1.19) 
  \psline[linewidth=0.04cm,arrowsize=0.1cm 2.0,arrowlength=1.4,arrowinset=0.4]{->}(5.2,1.23)(6.4,1.21) 
  \psline[linewidth=0.04cm,arrowsize=0.1cm 2.0,arrowlength=1.4,arrowinset=0.4]{->}(7.54,0.95)(8.3,0.47) 
  \psline[linewidth=0.04cm,arrowsize=0.1cm 2.0,arrowlength=1.4,arrowinset=0.4]{->}(5.24,-1.19)(6.42,-1.23) 
  \psline[linewidth=0.04cm,arrowsize=0.1cm 2.0,arrowlength=1.4,arrowinset=0.4]{->}(7.58,-0.93)(8.32,-0.53) 
  \psline[linewidth=0.04cm,arrowsize=0.1cm 2.0,arrowlength=1.4,arrowinset=0.4]{->}(2.78,-1.21)(4.0,-1.21) 
  \psline[linewidth=0.04cm,arrowsize=0.1cm 2.0,arrowlength=1.4,arrowinset=0.4]{->}(1.04,-0.45)(1.64,-0.93) 
  \psline[linewidth=0.04cm,linestyle=dashed,dash=0.16cm 0.16cm](2.1,0.55)(2.1,-0.55) 
  \psline[linewidth=0.04cm,linestyle=dashed,dash=0.16cm 0.16cm](4.62,0.57)(4.62,-0.55) 
  \psline[linewidth=0.04cm,linestyle=dashed,dash=0.16cm 0.16cm](6.98,0.59)(6.98,-0.59) 
  \usefont{T1}{ptm}{m}{n} 
  \rput(0.57703125,-0.08){0} 
  \usefont{T1}{ptm}{m}{n} 
  \rput(2.186875,1.18){1} 
  \usefont{T1}{ptm}{m}{n} 
  \rput(4.588594,1.18){$\cdots$} 
  \usefont{T1}{ptm}{m}{n} 
  \rput(6.9885936,1.2){\small{$m_0-1$}} 
  \usefont{T1}{ptm}{m}{n} 
  \rput(8.568594,-0.02){$m_0$} 
  \usefont{T1}{ptm}{m}{n} 
  \rput(4.608594,-1.22){$\cdots$} 
  \usefont{T1}{ptm}{m}{n} 
  \rput(2.1285937,-1.18){\tiny{$2m_0-1$}} 
  \usefont{T1}{ptm}{m}{n} 
  \rput(6.9685936,-1.24){\small{$m_0+1$}} 
  \end{pspicture}  
  }

      \caption{The path $\bi$. Vertices, that correspond to equal indices, are connected via dotted lines.} 
      \label{pict3} 
  
   \end{center} 
  
 \end{figure} 
 This path has $m_0+1$ distinct indices and it has $2m_0$ at all. It follows that there exist at least 2 one-element vertices of $\G_\bi$. Let
these one-element vertices be $\{s\}$ and $\{t\}$. Consider the pair  $(i_{t-1},i_{t})^\varepsilon{(t)}$. It must have an equal pair, but $i_t$ is not equal to any other index. It means, that $(i_{t-1},i_{t})^{\varepsilon(t-1)}=(i_{t},i_{t+1})^{\varepsilon(t)}$. It follows, that $\varepsilon(t-1) \neq \varepsilon(t)$. There are exactly two possibilities for this: $t=m_0$ or $t=0$. Assume without loss of generality that $s=0$ and $t=m_0$. Therefore $i_{m_0-1}=i_{m_0+1}$ (notice, that $(m_0-1)+(m_0+1)=2m_0$). Define $(m_0-1,1)$-path $\mathbf{j}$ as follows: $j_k:=i_k$ if $k \in \{0, \dots, m_0-2\}$, $j_{m_0-1}:=i_{m_0-1}=i_{m_0+1}$, $j_k:=i_{k+2}$ if $k \in \{m_0, \dots, 2(m_0-1)-1\}$. (See fig.\ref{pict4}) 
 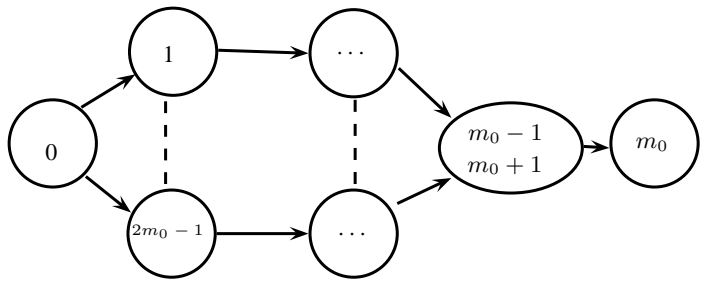
\begin{figure}[hbt] 
   \begin{center} 
  \scalebox{1} 
  { 
  \begin{pspicture}(0,-1.81)(9.2,1.81) 
  \pscircle[linewidth=0.04,dimen=outer](0.6,-0.01){0.6} 
  \pscircle[linewidth=0.04,dimen=outer](2.2,1.21){0.6} 
  \pscircle[linewidth=0.04,dimen=outer](4.6,1.19){0.6} 
  \pscircle[linewidth=0.04,dimen=outer](8.6,-0.01){0.6} 
  \pscircle[linewidth=0.04,dimen=outer](2.18,-1.21){0.6} 
  \pscircle[linewidth=0.04,dimen=outer](4.6,-1.21){0.6} 
  \psline[linewidth=0.04cm,arrowsize=0.1cm 2.0,arrowlength=1.4,arrowinset=0.4]{->}(0.98,0.47)(1.7,0.91) 
  \psline[linewidth=0.04cm,arrowsize=0.1cm 2.0,arrowlength=1.4,arrowinset=0.4]{->}(2.8,1.23)(4.0,1.19) 
  \psline[linewidth=0.04cm,arrowsize=0.1cm 2.0,arrowlength=1.4,arrowinset=0.4]{->}(2.78,-1.21)(4.0,-1.21) 
  \psline[linewidth=0.04cm,arrowsize=0.1cm 2.0,arrowlength=1.4,arrowinset=0.4]{->}(1.04,-0.45)(1.64,-0.93) 
  \psline[linewidth=0.04cm,linestyle=dashed,dash=0.16cm 0.16cm](2.1,0.55)(2.1,-0.55) 
  \psline[linewidth=0.04cm,linestyle=dashed,dash=0.16cm 0.16cm](4.62,0.57)(4.62,-0.55) 
  \usefont{T1}{ptm}{m}{n} 
  \rput(0.5840625,-0.12){0} 
  \usefont{T1}{ptm}{m}{n} 
  \rput(2.14375,1.18){1} 
  \usefont{T1}{ptm}{m}{n} 
  \rput(4.587188,1.18){$\cdots$} 
  \usefont{T1}{ptm}{m}{n} 
  \rput(8.567187,-0.02){$m_0$} 
  \usefont{T1}{ptm}{m}{n} 
  \rput(4.6071877,-1.22){$\cdots$} 
  \usefont{T1}{ptm}{m}{n} 
  \rput(2.1271875,-1.18){\tiny{$2m_0-1$}} 
  \psellipse[linewidth=0.04,dimen=outer](6.69,-0.07)(0.97,0.62) 
  \psline[linewidth=0.04cm,arrowsize=0.1cm 2.0,arrowlength=1.4,arrowinset=0.4]{->}(5.2,0.99)(5.9,0.33) 
  \psline[linewidth=0.04cm,arrowsize=0.1cm 2.0,arrowlength=1.4,arrowinset=0.4]{->}(5.18,-0.81)(5.9,-0.47) 
  \psline[linewidth=0.04cm,arrowsize=0.1cm 2.0,arrowlength=1.4,arrowinset=0.4]{->}(7.66,-0.05)(8.0,-0.07) 
  \usefont{T1}{ptm}{m}{n} 
  \rput(6.61625,0.12){$m_0-1$} 
  \usefont{T1}{ptm}{m}{n} 
  \rput(6.610781,-0.3){$m_0+1$} 
  \end{pspicture}  
  }

      \caption{The path $\mathbf{j}$  is a $(m_0-1,1)$-path} 
   \label{pict4} 
   \end{center} 
  
 \end{figure} 
 The path $\mathbf{j}$ is the $(m_0 - 1, 1)$-regular path. There is only one such path by inductive hypothesis and 
$(i_k=i_l)\Leftrightarrow(j_k=j_{l-2})\Leftrightarrow(k+(l-2)=2(m_0-1))\Leftrightarrow(k+l=2m_0)$.
 \end{proof} 
 
 \begin{defn} 
 Notice, that the vertex of a regular graph has two outgoing edges iff the corresponding index has the form  $i_{2mk}$ (because it should be
${(i_j,i_{j+1})^{\varepsilon(j)}=(i_j,i_{j+1})}$ and ${(i_{j-1},i_j)^{\varepsilon(j-1)}=(i_j,i_{j-1})}$ and it happens if and only if ${j=2mk}$). The distance between such vertex and vertex $\ver$  is called a \textit{type of vertex $\ver$}. The \textit{type of index} $i_j$ is the type of a vertex $\ver$ such that ${j\in\ver}$. It is clear, that index $i_j$ has type $j\pmod{2m}$ if $j\pmod{2m} \in \{ 0,\dots, m-1\}$ or type $-j\pmod{2m}$ in the other case. There are $m+1$ types of vertices. 
Note that only indices of the same type can be equal (this is proved in the case $p=1$ in Proposition \ref{3} and it will be proved for other cases below). 
 \label{def4} 
 \end{defn}

 Consider a collection of $(m,p_k)$-regular paths $(\bi_0,\bi_1,\dots,\bi_{m})$ (such that $\sum_{k=0}^{m}p_k=p-1$) and collection of corresponding
 regular graphs $(\G_0,\G_1,\dots,\G_{m})$.  Sum of path's lengths is $2m(p-1)$. We indicate the recipe how to
 obtain the $(m,p)$-regular graph from these collections. We take an $(m,1)$-regular graph and attach to its vertices the graphs from the collection in
 the following way: the graph $\G_0$ is attached to vertex of type $0$, the graph $\G_{1}$ is attached to vertex of type $1$, \dots, the
 graph $\G_{m}$ to the vertex of type $m$. 
For a more detailed argument we denote $\sum_{i=0}^k{p_i}$  by $P_k$ ($P_{m}=p-1$) and the indices of the
 $k^{th}$ path $\bi_{k}$ by $i^{(k)}_j$. The resulting $(m,p)$-regular graph is denoted by $\G_\mathbf{j}$. Define  the map 
$\Delta: \Delta(\bi_0,\bi_1,\dots,\bi_{m})=\mathbf{j}$ as follows 
 \begin{equation} 
 \begin{array}{ll} 
 j_0:=i_0^{(0)},j_1=i_1^{(0)}, \dots, j_{2mP_0-1}:=i^{(0)}_{2mp_0-1}, j_{2mP_0}:=j_0;\\ 
 j_{2mP_0+1}:=i_1^{(1)}, \dots, j_{2mP_1-1}:=i^{(1)}_{2mp_1-1}, j_{2mP_1}:=i^{(1)}_0, j_{2mP_1+1}:=j_{2mP_0+1};\\ 
 j_{2mP_1+2}:=i_2^{(2)},\dots,j_{2mP_2-1}:=i^{(2)}_{2mp_2-1},j_{2mP_2}:=i^{(2)}_0,j_{2mP_2+1}:=i^{(2)}_1,\\ 
 j_{2mP_2+2}:=j_{2mP_1+2};\\ 
 \dots\\ 
 j_{2mP_{m-1}+m}:=i_{m}^{(m)},\dots, j_{2m(p-1)-1}:=i^{(m)}_{2mp_{m}-1}, j_{2m(p-1)}:=i^{(m)}_{0},\\ 
 j_{2mp-m+1}:=j_{2mP_{m-1}+m-1},\dots, j_{2mp-k}:=j_{2mP_{k}+k}, \dots, j_{2mp-1}:=j_{2mP_{1}+1}. 
 \end{array} 
 \label{4} 
 \end{equation} 
 Let $\widetilde{\Delta}$ is the corresponding map $\widetilde{\Delta}(\G_0,\G_1,\dots,\G_{m})=\G_\mathbf{j}$.
 Graphically the constructon (\ref{4}) looks as follows: (fig. \ref{delta_pic}). 
 \begin{figure}[hbt] 
   \begin{center} 
  \scalebox{1} 
  { 
  \begin{pspicture}(0,-1.63)(11.52,1.65) 
  \psellipse[linewidth=0.04,dimen=outer](3.3,-0.97)(0.88,0.64) 
  \psellipse[linewidth=0.04,dimen=outer](5.66,-0.93)(0.88,0.64) 
  \psellipse[linewidth=0.04,dimen=outer](8.08,-0.95)(0.88,0.64) 
  \psarc[linewidth=0.04](3.29,0.62){0.95}{0.0}{180.0} 
  \psline[linewidth=0.04cm](2.34,0.63)(2.9,-0.39) 
  \psline[linewidth=0.04cm,arrowsize=0.2cm 2.0,arrowlength=1.4,arrowinset=0.4]{->}(4.24,0.63)(3.7,-0.45) 
  \psarc[linewidth=0.04](5.71,0.68){0.95}{0.0}{180.0} 
  \psline[linewidth=0.04cm](4.76,0.69)(5.32,-0.33) 
  \psline[linewidth=0.04cm,arrowsize=0.2cm 2.0,arrowlength=1.4,arrowinset=0.4]{->}(6.66,0.69)(6.12,-0.39) 
  \psarc[linewidth=0.04](8.11,0.64){0.95}{0.0}{180.0} 
  \psline[linewidth=0.04cm](7.16,0.65)(7.72,-0.37) 
  \psline[linewidth=0.04cm,arrowsize=0.2cm 2.0,arrowlength=1.4,arrowinset=0.4]{->}(9.06,0.65)(8.52,-0.43) 
  \psellipse[linewidth=0.04,dimen=outer](0.92,-0.99)(0.88,0.64) 
  \psarc[linewidth=0.04](0.95,0.6){0.95}{0.0}{180.0} 
  \psline[linewidth=0.04cm](0.0,0.61)(0.56,-0.41) 
  \psline[linewidth=0.04cm,arrowsize=0.2cm 2.0,arrowlength=1.4,arrowinset=0.4]{->}(1.9,0.61)(1.36,-0.47) 
  \psellipse[linewidth=0.04,dimen=outer](10.52,-0.95)(0.88,0.64) 
  \psarc[linewidth=0.04](10.55,0.64){0.95}{0.0}{180.0} 
  \psline[linewidth=0.04cm](9.6,0.65)(10.16,-0.37) 
  \psline[linewidth=0.04cm,arrowsize=0.2cm 2.0,arrowlength=1.4,arrowinset=0.4]{->}(11.5,0.65)(10.96,-0.43) 
  \psline[linewidth=0.04cm,arrowsize=0.2cm 2.0,arrowlength=1.4,arrowinset=0.4]{->}(1.84,-1.01)(2.44,-1.03) 
  \psline[linewidth=0.04cm,arrowsize=0.2cm 2.0,arrowlength=1.4,arrowinset=0.4]{->}(4.22,-0.99)(4.84,-1.01) 
  \psline[linewidth=0.04cm,arrowsize=0.2cm 2.0,arrowlength=1.4,arrowinset=0.4]{->}(6.58,-0.99)(7.24,-1.03) 
  \psline[linewidth=0.04cm,arrowsize=0.2cm 2.0,arrowlength=1.4,arrowinset=0.4]{->}(9.02,-0.99)(9.64,-1.03) 
  \usefont{T1}{ptm}{m}{n} 
  \rput(0.81171876,0.7){$\G_0$} 
  \usefont{T1}{ptm}{m}{n} 
  \rput(3.2717187,0.7){$\G_1$} 
  \usefont{T1}{ptm}{m}{n} 
  \rput(5.6517186,0.78){$\G_{\dots}$} 
  \usefont{T1}{ptm}{m}{n} 
  \rput(8.111719,0.8){$\G_{m-1}$} 
  \usefont{T1}{ptm}{m}{n} 
  \rput(10.491718,0.68){$\G_m$} 
  \usefont{T1}{ptm}{m}{n} 
  \rput(0.93703127,-0.84){0} 
  \rput(0.9340625,-1.26){$2mP_0$} 
   
  \rput(3.3140626,-0.84){\small{$2mP_0+1$}} 
  \rput(3.3140626,-1.2){\small{$2mP_1+1$}} 
  \usefont{T1}{ptm}{m}{n} 
  \rput(5.6340623,-1.0){$\cdots$} 
  \usefont{T1}{ptm}{m}{n} 
  \rput(8.17,-0.66){\tiny{$2mP_{m-1}+$}} 
  \rput(8.03,-0.86){\tiny{$+m-1$},} 
  \rput(8.07,-1.2){\tiny{$2mp-m+1$}} 
  \usefont{T1}{ptm}{m}{n} 
  \rput(10.55,-0.86){\footnotesize{$2mP_m-m$}} 
  \rput(10.55,-1.26){\footnotesize{$2mp-m$}} 
  \end{pspicture}} 
      \caption{The regular graph $\G_\mathbf{j}$, obtained from the collection of 
 $(m,p_i)$-regular graphs $(\G_0,\G_1,\dots,\G_{m})$} 
   \label{delta_pic} 
   \end{center} 
  
 \end{figure}
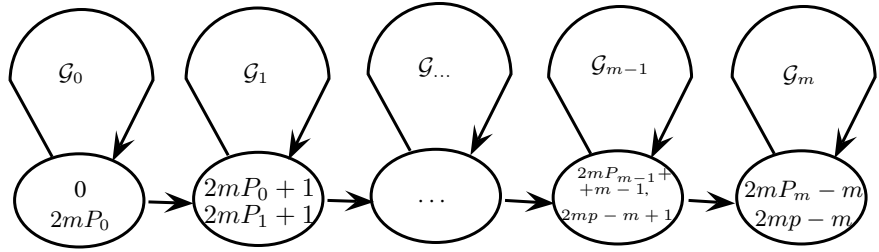 
  
 \textbf{Example.} 
 For example, we consider for $m=2$ the collection of $(2,p_k)$-regular graphs 
 $(\G_{2,2}$, $\G_{2,0}$, $\G_{2,1})$ (see fig. \ref{pict6}) and we obtain from it
a $(2,4)$-regular graph $\G_{2,4}$(see fig.\ref{pict7}). 
 \begin{figure}[ht] 
   \begin{center} 
    
  \scalebox{1} 
  { 
  \begin{pspicture}(0,-2.4)(9.04,2.42) 
  \pscircle[linewidth=0.04,dimen=middle](8.41,-1.79){0.59} 
  \pscircle[linewidth=0.04,dimen=middle](8.43,0.03){0.59} 
  \pscircle[linewidth=0.04,dimen=middle](8.43,1.81){0.59} 
  \psline[linewidth=0.04cm,arrowsize=0.2cm 2.0,arrowlength=1.4,arrowinset=0.4]{->}(8.42,-0.56)(8.44,-1.24) 
  \psline[linewidth=0.04cm,arrowsize=0.2cm 2.0,arrowlength=1.4,arrowinset=0.4]{->}(8.42,1.16)(8.44,0.58) 
  \pscircle[linewidth=0.04,dimen=outer](3.64,1.76){0.6} 
  \psellipse[linewidth=0.04,dimen=outer](3.6,-0.02)(1.8,0.6) 
  \pscircle[linewidth=0.04,dimen=outer](0.6,0.02){0.6} 
  \pscircle[linewidth=0.04,dimen=outer](6.62,0.02){0.6} 
  \pscircle[linewidth=0.04,dimen=outer](3.62,-1.8){0.6} 
  \psline[linewidth=0.04cm,arrowsize=0.2cm 2.0,arrowlength=1.4,arrowinset=0.4]{->}(3.62,-1.2)(3.64,-0.62) 
  \psline[linewidth=0.04cm,arrowsize=0.2cm 2.0,arrowlength=1.4,arrowinset=0.4]{->}(1.8,-0.02)(1.22,-0.02) 
  \psline[linewidth=0.04cm,arrowsize=0.2cm 2.0,arrowlength=1.4,arrowinset=0.4]{->}(3.62,1.12)(3.64,0.54) 
  \psline[linewidth=0.04cm,arrowsize=0.2cm 2.0,arrowlength=1.4,arrowinset=0.4]{->}(5.42,0.0)(6.04,-0.04) 
  \usefont{T1}{ptm}{m}{n} 
  \rput(3.591875,-0.07){1,3,5,7} 
  \usefont{T1}{ptm}{m}{n} 
  \rput(3.6240625,-1.87){0} 
  \usefont{T1}{ptm}{m}{n} 
  \rput(6.640625,-0.03){6} 
  \usefont{T1}{ptm}{m}{n} 
  \rput(3.631875,1.75){4} 
  \usefont{T1}{ptm}{m}{n} 
  \rput(0.6271875,-0.01){2} 
  \usefont{T1}{ptm}{m}{n} 
  \rput(8.424063,1.83){0} 
  \usefont{T1}{ptm}{m}{n} 
  \rput(8.406875,-0.03){1} 
  \usefont{T1}{ptm}{m}{n} 
  \rput(8.418593,-1.83){2} 
  \usefont{T1}{ptm}{m}{n} 
  \rput(1.0260937,-1.15){$\G_{2,2}$} 
  \usefont{T1}{ptm}{m}{n} 
  \rput(7.3915625,-1.73){$\G_{2,1}$} 
  \end{pspicture}  
  }		 
      \caption{Graphs $\G_{2,2}$ and $\G_{2,1}$. Graph $\G_{m,0}$ is empty} 
   \label{pict6} 
   \end{center} 
  
 \end{figure}
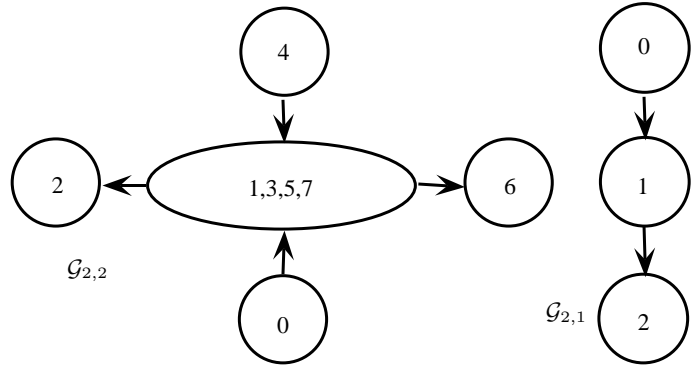 
 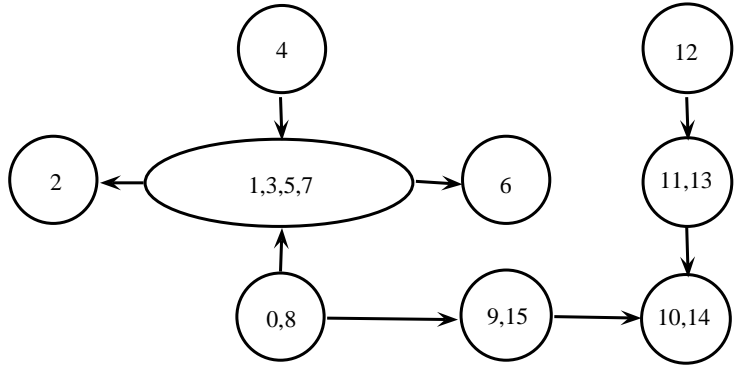
\begin{figure}[ht] 
   \begin{center} 
    
   	  \scalebox{1} 
  		{ 
  		\begin{pspicture}(0,-2.41)(9.62,2.41) 
  		\pscircle[linewidth=0.04,dimen=middle](3.62,1.79){0.58} 
  		\psellipse[linewidth=0.04,dimen=middle](3.58,0.01)(1.78,0.58) 
  		\pscircle[linewidth=0.04,dimen=middle](0.58,0.05){0.58} 
  		\pscircle[linewidth=0.04,dimen=middle](6.6,0.05){0.58} 
  		\pscircle[linewidth=0.04,dimen=middle](3.6,-1.77){0.58} 
  		\psline[linewidth=0.04cm,arrowsize=0.2,arrowlength=1.4,arrowinset=0.4]{->}(3.6,-1.17)(3.62,-0.59) 
  		\psline[linewidth=0.04cm,arrowsize=0.2,arrowlength=1.4,arrowinset=0.4]{->}(1.78,0.01)(1.2,0.01) 
  		\psline[linewidth=0.04cm,arrowsize=0.2,arrowlength=1.4,arrowinset=0.4]{->}(3.6,1.15)(3.62,0.57) 
  		\psline[linewidth=0.04cm,arrowsize=0.2,arrowlength=1.4,arrowinset=0.4]{->}(5.4,0.03)(6.02,-0.01) 
  		\usefont{T1}{ptm}{m}{n} 
  		\rput(3.6009376,-0.04){1,3,5,7} 
  		\usefont{T1}{ptm}{m}{n} 
  		\rput(6.5953126,-0.02){6} 
  		\usefont{T1}{ptm}{m}{n} 
  		\rput(3.6209376,1.78){4} 
  		\usefont{T1}{ptm}{m}{n} 
  		\rput(0.61859375,0.02){2} 
  		\pscircle[linewidth=0.04,dimen=middle](6.6,-1.75){0.6} 
  		\usefont{T1}{ptm}{m}{n} 
  		\rput(3.613125,-1.82){0,8} 
  		\pscircle[linewidth=0.04,dimen=middle](8.99,-1.8){0.59} 
  		\pscircle[linewidth=0.04,dimen=middle](9.01,0.02){0.59} 
  		\pscircle[linewidth=0.04,dimen=middle](9.01,1.8){0.59} 
  		\psline[linewidth=0.04cm,arrowsize=0.2,arrowlength=1.4,arrowinset=0.4]{->}(4.22,-1.79)(5.96,-1.81) 
  		\psline[linewidth=0.04cm,arrowsize=0.2,arrowlength=1.4,arrowinset=0.4]{->}(7.24,-1.77)(8.4,-1.79) 
  		\psline[linewidth=0.04cm,arrowsize=0.2,arrowlength=1.4,arrowinset=0.4]{->}(9.0,-0.57)(9.02,-1.25) 
  		\psline[linewidth=0.04cm,arrowsize=0.2,arrowlength=1.4,arrowinset=0.4]{->}(9.0,1.15)(9.02,0.57) 
  		\usefont{T1}{ptm}{m}{n} 
  		\rput(6.6209373,-1.78){9,15} 
  		\usefont{T1}{ptm}{m}{n} 
  		\rput(8.962656,-1.8){10,14} 
  		\usefont{T1}{ptm}{m}{n} 
  		\rput(8.993906,0.02){11,13} 
  		\usefont{T1}{ptm}{m}{n} 
  		\rput(8.991406,1.76){12} 
  		\end{pspicture}  
  		}

      \caption{The resulting graph $\G_{2,4}$} 
   \label{pict7} 
   \end{center} 
  
 \end{figure} 
  
  \begin{prop} 
Using the above construction we get an  $(m,p)$-regular graph. 
  \end{prop} 
  
 \begin{proof} 
 Indeed, the graph $\G_\mathbf{j}$ has exactly $mp$ edges and $k_r=2$ for all $r=1,2,\dots,mp$. 
Furthermore, there are exactly $mp+1$ vertices (there is
no new-introduced vertex and there are exactly $\sum_{i=0}^{m}{(mp_i+1)}=m(p-1)+m+1=mp+1$ vertices of graphs $\G_k$). 
 \end{proof} 
  
Note, that the map ${\Delta}$ is  the injection. 

 Now we consider the arbitrary $(m,p)$-regular path  $\bi$  and try to construct inverse map for $\Delta$. Denote 
 \begin{equation} 
 \begin{array}{ll} 
 	J_0:=\{ j:i_j=i_0 \}; \\ 
 	J_k:=\{ j: j\neq 2mp-k, i_j=i_{2mp-k}\} ,~ k \in \{ 1, \dots, m-1 \} ;\\ 
 	J_{m}:=\{ j:i_j=i_{2mp-m} \};\\ 
 	\overline{J_k}:=\max(J_k),~ \underline{J_k}:=\min(J_k). 
 \end{array} 
 \label{5} 
 \end{equation} 
 We will prove that the sets $J_k$ have some remarkable properties and after that it will be clear, how to obtain a collection of regular paths from one (big) regular path. 
 \begin{prop} 
  $J_k$ is nonempty. 
 \label{6} 
 \end{prop} 
  \begin{proof} 
 Indeed, there is $0 \in J_0$ and $2mp-m \in J_{m}$.
 If $J_k$ is void with ${k \in \{ 1,\dots,m-1 \}}$,  then the index $i_{2mp-k}$
 has no equal indices  in the path $\bi$. But in this case the pair 
$(i_{2mp-k-1},i_{2mp-k})^-$ has no equivalent for the following reason.
 The index $i_{2mp-k}$
 appears in $(2mp-k,2mp-k+1)^-$  and in $(2mp-k-1,2mp-k)^-$ only and they are not equivalent. 
But each edge in a regular path has equivalent one, a contradiction. 
Therefore the initial assumption that $J_k$ is void must be false.
  \end{proof} 
  
 \begin{prop} 
   $J_k$ ($0 \leq k \leq m$) are pairwise disjoint and if $k<l$ then $J_k < J_l$ (for all $j\in J_k$ and  for all $i \in J_l$ the inequality $j<i$
holds). 
 \label{7} 
 \end{prop} 
   
 \begin{proof} 
 Indeed, if $J_k \bigcap J_l \neq \emptyset$ then $i_{2mp-k}=i_{2mp-l}$.
The  edges of the path $\bi$ have the same orientation on the section 
$(2mp-k,2mp-k-1)^-, \dots, (2mp-l+1,2mp-l)$ and therefore 
the graph $G_\bi$ has a cycle. But a regular graph is a tree, a contradiction.
 Thus $J_k \bigcap J_l = \emptyset$. We prove the second part of Proposition \ref{7} 
for the case $l=k+1$ only (which is sufficient). Consider the edge $(2mp-(k+1),2mp-k)^-$. 
It must be equivalent to an edge $(t,t+1)^+$ with some $t \in J_k$ and $t+1 \in J_{k+1}$.
 If there exists $s \in J_k$ such that $s>t$ then $s>t+1$ ($J_k \bigcap J_{k+1} = \emptyset$).
 The edge $(t,t+1)$ is not equivalent to any edge in the section ${(t+1,t+2),\dots,(s-1,s)}$,
 because it has only one equivalent edge $(2mp-(k+1),2mp-k)^-$. 
It follows that there are two different 
 paths in the graph $G_\bi$ which 
 connect vertex $\mathcal{U}$ (such that $t+1\in\mathcal{U}$)  and vertex $\ver$ (such that $s\in\ver$ and $t\in\ver$), that is there is a cycle in the 
the graph $G_\bi$, and hence there is a contradiction. 
Therefore $t=\max J_k= \overline{J_k}$. Similarly, $t+1=\underline{J_{k+1}}$. 
 It follows, that  $\overline{J_k}+1=\underline{J_{k+1}}$ and for all$j \in J_k$ 
and for all $i \in J_{k+1}$ the inequality $j<i$ holds. 
 \end{proof} 
 \begin{prop} 
 For all $k$ the difference $(\overline{J_k}-\underline{J_k})$ is divisible by $2m$.
 \label{8} 
 \end{prop} 
  
 \begin{proof} 
 Denote $(\overline{J_k}-\underline{J_k})\pmod{2m}$ by $d_k$. Notice, that  $(\overline{J_k}-\underline{J_k})$ is the number of edges in the path's 
section   ${(\underline{J_k},\underline{J_k}+1),\dots, (\overline{J_k}-1,\overline{J_k})}$. Notice, that the orientation of edges changes after every $m$ steps. Edges of the form $(i_{\overline{J_k}},i_{\underline{J_{k+1}}})^+$ have the same orientation. It follows, that $d_0 \leq m-1$, $d_1 \leq m-1$, $(d_0+1+d_1) \leq m-1~(mod~2m)$ (and so $(d_0+1+d_1) \leq m-1$, because $0 \leq d_0+1+d_1 \leq 2m-1$), $\dots$, $0 \leq d_0+1+d_1+1+...+d_{m-2}+1+d_{m-1}\leq m-1$ (similarly), i.e. $0 \leq \sum_{k=0}^{m-1}{d_k}+m-1 \leq m-1$. Therefore, $d_k=0$ for all ${k= 0,1,\dots,m-1}$. Consider all edges of the path $\bi$.
 There are $m$ edges of the form $(\overline{J_k},\underline{J_{k+1}})$, 
$m$ edges of the form $(2mp-k,2mp-k+1)^-$ with some ${k=1,2,\dots,m}$ and 
all the remaing ones  are in sections of 
the form ${(\underline{J_k},\underline{J_k}+1),\dots, (\overline{J_k}-1,\overline{J_k})}$.
 There are $2mp$ edges in total. 
Therefore, $\sum_{k=0}^{m}(\overline{J_k}-\underline{J_k})+m+m=2mp$ 
and hence $\sum_{k=0}^{m} d_k = 0\pmod{2m}$. It follows, that $d_{m}=0$ too.
 \end{proof} 
 \begin{prop} 
 If $\underline{J_k}<t<\overline{J_k}$ and $\underline{J_l}<s<\overline{J_l}$, then $i_t \neq i_s$. In other words, sections of the path $\bi$ of the form $(\underline{J_k},\underline{J_k}+1), \dots, (\overline{J_k}-1,\overline{J_k})$ with $k=0,1,\dots,m$ are disjoint.
 \label{10} 
 \end{prop} 
 \begin{proof} 
 Without loss of generality we consider $l>k$. Assume that $i_t=i_s$. 
In this case the section $(\underline{J_k},\underline{J_k}+1), \dots, (s-1,s)$ 
contains  the edge $(\overline{J_k},\underline{J_{k+1}})$, 
and the section $(t,t+1), \dots, (\overline{J_k}-1,\overline{J_k})$ does not contain it
 or its equivalent $(2mp-k-1,2mp-k)^-$. Thus,  there are two non-equal paths in the regular 
graph $G_\bi$ which  connected vertex $\mathcal{U}$ (such that $\overline{J_k}\in\mathcal{U}$)  and vertex $\ver$ (such that $s\in\ver$ and $t\in\ver$), that is there is a cycle in the 
the graph $G_\bi$. Therefore, the initial assumption must be false. 
 \end{proof} 
  
 Now we can describe the inverse map for $\Delta$. Let 
$p_k:= (\overline{J_k}-\underline{J_k})/2m$ ($p_k$ is a nonnegative integer by Proposition 
\ref{8}). Furthermore, we have for  sum $\sum_{k=0}^{m-1}{p_k}=p-1$ 
(see the proof of Proposition \ref{8}).
 Denote by $\mathbf{j}^{(k)}$ the $k$-th resulting path 
(it has a length $2mp_k$  and if $p_k=0$ then $\mathbf{j}_k$ is empty). Let  
  
 \begin{equation} 
 \label{9} 
 j^{(k)}_{t}:=i_{\underline{J_k}+((t-k)~mod~2mp_k))},\, t \in \{0,\dots,2mp_k-1\},\, k \in \{0,\dots,m\}. 
 \end{equation} 
 Now one obtains the collection $(\G_{2,2},~\G_{2,0},~ G_{2,1})$ (see fig. \ref{pict6}) from the graph $\G_{2,4}$ (see fig.\ref{pict7}) in the way described in (\ref{9}). 
 \begin{prop} 
 The collection of paths $(\mathbf{j}^{(0)},\mathbf{j}^{(1)},\dots,\mathbf{j}^{(m)})$ 
$($defined by $(\ref{9}))$ is the collection of regular paths.   
 \end{prop} 
  
 \begin{proof} 
 In fact, the path $\mathbf{j}^{(k)}$ is almost the same as the section  
$(\underline{J_k},\underline{J_k}+1),\dots,(\overline{J_k}-1,\overline{J_k})$ 
of the regular path $\bi$ . This section contains $2mp_k$ edges. 
Each of these edges has an equivalent one in the same section by Proposition \ref{10}.
 Therefore this section contains exactly $mp_k+1$ distinct indices because of connectivity 
and  non-cyclicity. Hence the path $\mathbf{j}^{(k)}$ is a regular path. 
  \end{proof} 
 Thus, $\widetilde{\Delta}$ is the bijection between  $T^{reg}_{m,p}$ and 
$\bigcup T^{reg}_{m,p_0}\times T^{reg}_{m,p_1}\times\cdots\times T^{reg}_{m,p_m}$, 
where the union is taken over all $p_0+p_1+\dots+p_m=p-1$. 
Hence $\#T_{m,p}^{reg}=M_p^{(m)}$ and Lemma \ref{lem2} is proved.
 \end{proof} 
 Lemmas \ref{lem1} and \ref{lem2} show that the moments of the spectral distribution converge to $M^{(m)}_p$. Thus Theorem \ref{main_thm} is
proved.


\begin{thebibliography}{99} 
\bibitem{Bai99} Z. D. Bai, \emph{Methodologies in spectral analysis of large-dimensional random
matrices, a review,} Statist. Sinica 9 (1999), no. 3, 611-677.
\bibitem{banica08} Banica, T. Belinschi, S. Capitaine,   M. and Collins B. 
{\em Free Bessel Laws}
Preprint. arXiv:0710.5931
\bibitem{Carl1926} Carleman T. \emph{Les fonctions quasi-analytiques}, Paris, 1926.
\bibitem{Ges2006} Gessel, Ira M.; Xin, Guoce \emph{The generating function of ternary trees and continued fractions.}  Electron. J. Combin.  13  (2006),  no. 1.
\bibitem{CMath} Ronald L. Graham, Donald E. Knuth, Oren Patashnik, \emph{Concrete Mathematics: A Foundation for Computer Science} 
\bibitem{MarPas}  Marchenko and V., Pastur, L.
{\em The eigenvalue distribution in some ensembles of random
matrices}.\newline Math.USSR Sbornik, {\bf1} (1967), 457-483
\bibitem{Ml09} W. Mlotkowski, \emph{Fuss--Catalan numbers in noncommutative probability}, preprint. 
\bibitem{NSp} A. Nica, R. Speicher, \emph{Lectures on the Combinatorics of Free Probability}, Cambridge University Press, 2006 
 
\bibitem{ora} F. Oravecz, \emph{On the powers of Voiculescu�s circular element}, Studia Math. 145 (2001) 
\bibitem{Wig58} Wigner, E.
     {\em On the distribution of the roots of certain symmetric matrices}.\newline
     Ann. of Math.{\bf67} (1958), 325--327.
\end{thebibliography}
\end{document}